\documentclass[12pt]{article}
\usepackage{latexsym}
\usepackage{amsmath}
\usepackage{amssymb}
\usepackage{amsfonts}
\usepackage[colorlinks=true]{hyperref}
\hypersetup{urlcolor=blue, citecolor=red}

\numberwithin{equation}{section}

\newtheorem{theorem}{Theorem}[section]
\newtheorem{lemma}[theorem]{Lemma}

\newtheorem{remark}[theorem]{Remark}
\newcommand{\eproof}{{\mbox{\ }~\hfill
\mbox{\large $\Box$} \par \vskip 10pt}}

\newcommand{\supp}{\mbox{\rm supp\,}}

\newcommand{\R}{{\mathbb R}}

\newcommand{\pf}{\noindent{\bf Proof}}

\newcommand{\curl}{{\rm curl}}

\title{Asymptotic behavior of solutions of the stationary Navier-Stokes equations in an exterior domain}
\author{Ching-Lung Lin\thanks{Department of Mathematics, NCTS, National Cheng Kung University,
Tainan 701, Taiwan. Email:cllin2@mail.ncku.edu.tw}\quad Gunther Uhlmann\thanks{Department of Mathematics, University of Washington, Seattle, WA 98195-4350, USA
and Department of Mathematics University of California, Irvine, CA 92697-3875, USA. Email:gunther@math.washington.edu}\quad
Jenn-Nan Wang\thanks{Department of Mathematics, Taida Institute of
Mathematical Sciences, NCTS (Taipei), National Taiwan University,
Taipei 106, Taiwan. Email: jnwang@math.ntu.edu.tw}}

\date{}

\begin{document}
\maketitle

\begin{abstract}
We study the asymptotic behavior of an incompressible fluid around a
bounded obstacle. The problem is modeled by the stationary
Navier-Stokes equations in an exterior domain in $\R^n$ with $n\ge
2$. We will show that, under some assumptions, any nontrivial
velocity field obeys a minimal decaying rate $\exp(-Ct^2\log t)$ at
infinity. Our proof is based on appropriate Carleman estimates.
\end{abstract}

\section{Introduction}\label{sec1}
\setcounter{equation}{0}

Let $B$ be a bounded domain in $\R^n$ and $\Omega=\R^n\setminus\bar{B}$ with $n\ge 2$. Without loss of generality, we assume $0$ is in the interior of $B$ and
$B\subset B_1(0)=\{x: |x|<1\}$. Assume that $\Omega$ is filled with an incompressible fluid described by the stationary Navier-Stokes equations
\begin{equation}\label{1.1}
\begin{cases}
\begin{array}{l}
-\Delta u+u\cdot\nabla u+\nabla p=0\quad\text{in}\quad\Omega,\\
\nabla\cdot u=0\quad\text{in}\quad\Omega.
\end{array}
\end{cases}
\end{equation}
We are interested in the following question: under some boundedness assumption on $u$, what is the minimal decaying rate at infinity of any nontrivial $u$ satisfying \eqref{1.1}?

To put the problem in perspective, we first mention some related
results. In three dimensions, Finn \cite{finn1} showed that if
$u|_{\partial B}=0$ and $u=o(|x|^{-1})$, then $u$ is trivial. On the
other hand,  Dyer and Edmunds \cite{dyed} proved that if $u$ is
$C^2$ bounded and $u=O(\exp(-\exp(\alpha|x|^3)))$ for all
$\alpha>0$, then $u$ is trivial. We remark that in Finn's result $u$
is required to satisfy the homogeneous Dirichlet condition on
$\partial B$ and a decaying condition at infinity, while in Dyer and
Edmunds's result, with the assumption of $C^2$ boundedness, only the
local behavior of $u$ at infinity is needed. We have showed in an
early paper \cite{luw1} that for $n=2$ or $3$, if $u$ is bounded in
$\Omega$, then any nontrivial $u$ of \eqref{1.1} can not decay
faster than certain double exponential at infinity (see
\cite[Corollary~1.6]{luw1} for details). In the present paper, we
improve significantly on that result, and the result of \cite{dyed}
by showing that the minimal decaying rate of any nontrivial $u$ is
close to exponential in dimension $n\ge 2$. We now state the main
theorem of the paper. Denote
\begin{equation*}
M(t)=\inf_{|x|=t}\int_{|y-x|<1}|u(y)|^2 dy.
\end{equation*}
\begin{theorem}\label{thm1.1}
Let $u\in (H^1_{loc}(\Omega))^{n}$ be a nontrivial solution of \eqref{1.1} with an appropriate $p\in H^1_{loc}(\Omega)$. Assume that
\begin{equation}\label{2bound}
\|u\|_{L^{\infty}(\Omega)}\leq \lambda\quad\text{if}\quad n=2,
\end{equation}
or
\begin{equation}\label{3bound}
\|u\|_{L^{\infty}(\Omega)}+\|\nabla u\|_{L^{\infty}(\Omega)}\leq \lambda\quad\text{if}\quad n\ge 3.
\end{equation}
Then there exist $C>0$ depending on $\lambda$, $n$, and $\tilde t>0$ depending on $\lambda$, $n$, $M(10)$ such that
\begin{equation*}
M(t)\geq \exp(-Ct^2\log t)\quad\text{for}\quad t\ge \tilde t.
\end{equation*}
\end{theorem}
\begin{remark}
It is interesting to compare our result with a similar result for the Schr\"odinger equation proved by Bourgain and Kenig \cite{bou} (see also \cite{kenigcna}). In \cite{bou}, they considered the Schr\"odinger equation
$$\Delta u+V(x)u=0\quad\text{in}\quad\R^n.$$ Under the assumption that $|V|\le 1$, $|u_0|\le C_0$, and $u(0)=1$, they proved
$$\inf_{|x_0|=R}\sup_{B(x_0,1)}|u(x)|\ge C\exp(-R^{4/3}\log R)\quad\text{for}\quad R>>1.$$
\end{remark}

We prove our results by using appropriate Carleman estimates. We will
use weights which are slightly less singular than negative
powers of $|x|$ (see estimates \eqref{2.1}). The method of obtaining a decaying rate is a detour from that of deriving three-ball inequalities using Carleman estimates.

This paper is organized as follows. In Section~2, we reduce the Navier-Stokes equations to a new system by using the vorticity function. Then we state some suitable
Carleman-type estimates. A technical interior estimate is proved in
Section~3. Section~4 is devoted to the proof of Theorem~\ref{thm1.1}.

\section{Reduced system and Carleman estimates}\label{sec2}
\setcounter{equation}{0}

Fixing $x_0$ with  $|x_0|=t>>1$, we define
\begin{equation*}
\begin{array}{l}
w(x)=(at)u(atx+x_0),\quad\tilde{p}(x)=(at)^2p(atx+x_0),
\end{array}
\end{equation*}
where $a=8/s$ and $0<s<8$ is a small constant which will be determined in the proof of Theorem \ref{thm1.1}. Likewise, we denote
$$\Omega_{t}:=B_{\frac{1}{a}-\frac{1}{at}}(0).$$  From \eqref{1.1}, it is easy to check that
\begin{equation}\label{2.1}
\begin{cases}
\begin{array}{l}
-\Delta w+w\cdot\nabla w+\nabla \tilde{p}=0\quad\text{in}\quad\Omega_{t},\\
\nabla\cdot w=0\quad\text{in}\quad\Omega_{t}.
\end{array}
\end{cases}
\end{equation}
In view of \eqref{2bound} and \eqref{3bound}, we have that
\begin{equation}\label{2.2}
\|w\|_{L^{\infty}(\Omega_{t})}\leq at\lambda
\end{equation}
or
\begin{equation}\label{2.3}
\|w\|_{L^{\infty}(\Omega_{t})}\leq at\lambda\quad\text{and}\quad \|\nabla w\|_{L^{\infty}(\Omega_{t})}\le(at)^2\lambda.
\end{equation}

To study the Navier-Stokes equation, it is often advantageous to consider the vorticity equation. Let us now define the vorticity $q$ of the velocity $w$ by
$$q=\curl w:=\frac{1}{\sqrt{2}}(\partial_iw_j-\partial_jw_i)_{1\le i,j\le n}.$$ Note that here $q$ is a matrix-valued function. The formal transpose of $\curl$ is given by
$$(\curl^{\top}v)_{1\le i\le n}:=\frac{1}{\sqrt{2}}\sum_{1\le j\le n}\partial_j(v_{ij}-v_{ji}),$$ where $v=(v_{ij})_{1\le i,j\le n}$. It is easy to see that
\begin{equation*}
\Delta w=\nabla(\nabla\cdot w)-\curl^{\top}\curl w
\end{equation*}
(see, for example, \cite{mm} for a proof), which implies
\begin{equation*}
\Delta w+\curl^{\top}q=0\quad\text{in}\quad\Omega_{t}.
\end{equation*}
On the other hand, we observe that
$$w\cdot\nabla w=\nabla(\frac 12 |w|^2)-\sqrt{2}(\curl w)w=\nabla(\frac 12|w|^2)-\sqrt{2}qw.$$ Thus, applying $\curl$ on the first equation of \eqref{2.2}, we have that
\begin{equation*}
-\Delta q+Q(q)w+q(\nabla w)^{\top}-(\nabla w)q^{\top}=0\quad\text{in}\quad\Omega_{t},
\end{equation*}
where $$(Q(q)w)_{ij}=\sum_{1\le k\le n}(\partial_jq_{ik}-\partial_iq_{jk})w_k.$$ Now for $n=2$, due to $\nabla\cdot w=0$, it is easily seen that
$$q(\nabla w)^{\top}-(\nabla w)q^{\top}=0.$$
Therefore, we will consider the system
\begin{equation}\label{2.4}
\begin{cases}
\begin{array}{l}
-\Delta q+Q(q)w+q(\nabla w)^{\top}-(\nabla w)q^{\top}=0\quad\text{in}\quad\Omega_{t},\\
\Delta w+\curl^{\top} q=0\quad\text{in}\quad\Omega_{t}
\end{array}
\end{cases}
\end{equation}
for $n\ge 3$, and
\begin{equation}\label{2.5}
\begin{cases}
\begin{array}{l}
-\Delta q+Q(q)w=0\quad\text{in}\quad\Omega_{t},\\
\Delta w+\curl^{\top} q=0\quad\text{in}\quad\Omega_{t}
\end{array}
\end{cases}
\end{equation}
for $n=2$. In order to prove the main theorem, putting together \eqref{2.4}, \eqref{2.5}, and using \eqref{2.2}, \eqref{2.3}, it suffices to consider
\begin{equation}\label{2.6}
\begin{cases}
\begin{array}{l}
-\Delta q+A(x)\cdot\nabla q+B(x)q=0\quad\text{in}\quad\Omega_{t},\\
\Delta w+\curl^{\top} q=0\quad\text{in}\quad\Omega_{t}
\end{array}
\end{cases}
\end{equation}
with
$$\|A\|_{L^{\infty}(\Omega_t)}\le at\lambda\quad\text{and}\quad\|B\|_{L^{\infty}(\Omega_t)}\le (at)^2\lambda.$$

For our proof, we will apply Carleman estimates with weights $\varphi_{\beta}=\varphi_{\beta}(x) =\exp
(-\beta\tilde{\psi}(x))$, where $\beta>0$ and $\tilde{\psi}(x)=\log
|x|+\log((\log |x|)^2)$.
\begin{lemma}\label{lem2.1}
There exist a sufficiently small number $r_1>0$,
a sufficiently large number $\beta_1>2$, a positive constant $C$, such that
for all $v\in U_{r_1}$ and $\beta\geq \beta_1$, we have that
\begin{equation}\label{2.7}
\begin{array}{l}
\int\varphi^2_\beta (\log|x|)^2(\beta|x|^{4-n}|\nabla v|^2+\beta^3|x|^{2-n}|v|^2)dx
\le C\int \varphi^2_\beta(\log|x|)^{4}|x|^{6-n}|\Delta v|^2dx,
\end{array}
\end{equation}
where $U_{r_1}=\{v\in C_0^{\infty}(\R^n\setminus\{0\}): \mbox{\rm
supp}(v)\subset B_{r_1}\}$.
\end{lemma}

Lemma~\ref{lem2.1} is a modified form of \cite[Lemma~2.4]{lin3}. For the sake of brevity, we omit the proof here. Applying Lemma \ref{lem2.1} with $\beta=\beta+1$, we have the following Carleman estimates.
\begin{lemma}\label{lem2.2}
There exist a sufficiently small number $r_1>0$,
a sufficiently large number $\beta_1>1$, a positive constant $C$, such that
for all $v\in U_{r_1}$ and $\beta\geq \beta_1$, we have
\begin{equation}\label{2.8}
\int \varphi^2_\beta (\log|x|)^{-2}|x|^{-n}(\beta |x|^{2}|\nabla
v|^2+\beta^3|v|^2)dx\le C \int \varphi^2_\beta|x|^{-n}(|x|^{4}|\Delta v|^2)dx.
\end{equation}
\end{lemma}

\section{Interior estimates}\label{sec3}
\setcounter{equation}{0}

In addition to Carleman estimates, we also need the following
interior inequality.

\begin{lemma}\label{lem3.1}
For any $0<a_1<a_2$ such that $B_{a_2}\subset\Omega_{t}$, let $X=B_{a_2}\backslash \bar{B}_{a_1}$ and $d(x)$ be
the distance from $x\in X$ to $\mathbb{R}^n\backslash X$. We have
\begin{eqnarray}\label{3.1}
&&\int_{X}d(x)^{2}|\nabla w|^2dx+\int_{X}d(x)^{4}|\nabla q|^2dx+\int_{X}d(x)^{2}|q|^2dx\notag\\
&\le&C(at )^{12}\int_{X}|w|^2 dx.
\end{eqnarray}
where the constant $C$ is independent of $r$, $a$, $t$ and $(w,q)$.
\end{lemma}

\pf.
By
elliptic regularity, we obtain from \eqref{1.1} that $u\in
H^{2}_{loc}(\Omega)$ and hence $w\in H^2(\Omega_t)$. It is trivial that
\begin{equation}\label{3.2}
\begin{array}{l}
\| v\|_{H^{1}(\mathbb{R}^n)}\lesssim \|\Delta
v\|_{L^{2}(\mathbb{R}^n)} +\|v\|_{L^{2}(\mathbb{R}^n)}
\end{array}
\end{equation}
for all $v\in H^{2}(\mathbb{R}^n)$. By changing variables $x\to
E^{-1}x$ in \eqref{3.2}, we obtain

\begin{equation}\label{3.3}
\begin{array}{l}
\sum_{|\alpha|\leq 1}E^{2-|\alpha|}\|
D^{\alpha}v\|_{L^{2}(\mathbb{R}^n)}\lesssim(\|\Delta
v\|_{L^{2}(\mathbb{R}^n)} +E^2\|v\|_{L^{2}(\mathbb{R}^n)})
\end{array}
\end{equation}
for all $v\in H^{2}(\mathbb{R}^n)$. To apply \eqref{3.3} to $w$, we
need to cut-off $w$. So let $\xi(x)\in C^{\infty}_0 ({\mathbb R}^n)$
satisfy $0\le\xi(x)\leq 1$ and
\begin{equation*}
\xi (x)=
\begin{cases}
\begin{array}{l}
1,\quad |x|<1/4,\\
0,\quad |x|\geq 1/2.
\end{array}
\end{cases}
\end{equation*}
Let us denote $\xi_y(x)=\xi((x-y)/d(y))$.   For $y\in X$, we apply
\eqref{3.3} to $\xi_y(x) w(x)$  and use the second equation of \eqref{2.6} to get
that
\begin{eqnarray}\label{3.4}
&& E^{2}\int_{|x-y|\leq d(y)/4}|\nabla w|^2dx\notag\\
&\le&C_1\int_{|x-y|\leq d(y)/2}|\nabla q|^2dx+C_1\int_{|x-y|\leq d(y)/2}d(y)^{-2}|\nabla w|^2dx\notag\\
&& +C_1(E^4+d(y)^{-4})\int_{|x-y|\leq d(y)/2}|w|^2dx.
\end{eqnarray}
Now taking $E=M^3d(y)^{-1}$ for some constant $M>1$ and
multiplying $d(y)^{4}$ on both sides of \eqref{3.4}, we have
\begin{eqnarray}\label{3.5}
&&M^{6}d(y)^{2}\int_{|x-y|\leq d(y)/4}|\nabla w|^2dx\notag\\
&\le& C_1\int_{|x-y|\leq d(y)/2}d(y)^{4}|\nabla
q|^2dx+C_1\int_{|x-y|\leq
d(y)/2}d(y)^{2}|\nabla w|^2dx\notag\\
&&+C_1(M^{12}+1)\int_{|x-y|\leq d(y)/2}|w|^2dx.
\end{eqnarray}

Integrating $d(y)^{-n}dy$ over $X$ on both sides of \eqref{3.5} and
using Fubini's Theorem, we get that
\begin{eqnarray}\label{3.6}
&&M^{6}\int_{X}\int_{|x-y|\leq d(y)/4}d(y)^{2-n}|\nabla w|^2dydx\notag\\
&\le&C_1\int_{X}\int_{|x-y|\leq d(y)/2}d(y)^{4-n}|\nabla q(x)|^2dydx\notag\\
&&\quad+ C_1\int_{X}\int_{|x-y|\leq d(y)/2}d(y)^{2-n}|\nabla w|^2dydx\notag\\
&&\quad+2C_1M^{12}\int_{X}\int_{|x-y|\leq d(y)/2}d(y)^{-n}|w|^2dydx.
\end{eqnarray}
Note that $|d(x)-d(y)|\leq |x-y|$. If $ |x-y|\leq d(x)/3$, then
\begin{equation}\label{3.7}
\begin{array}{l}
2d(x)/3\leq d(y)\leq 4d(x)/3.
\end{array}
\end{equation}
On the other hand, if $ |x-y|\leq d(y)/2$, then
\begin{equation}\label{3.8}
\begin{array}{l}
d(x)/2\leq d(y)\leq 3d(x)/2.
\end{array}
\end{equation}
By \eqref{3.7} and  \eqref{3.8}, we have
\begin{equation}\label{3.9}
\begin{cases}
\begin{array}{l}
\int_{|x-y|\leq d(y)/4}d(y)^{-n}dy\geq (3/4)^n\int_{|x-y|\leq
d(x)/6}d(x)^{-n}dy\geq 8^{-n}\int_{|y|\leq 1}dy,\\
\int_{|x-y|\leq d(y)/2}d(y)^{-n}dy\leq 2^n\int_{|x-y|\leq
3d(x)/4}d(x)^{-n}dy\leq (3/2)^{n}\int_{|y|\leq 1}dy.
\end{array}
\end{cases}
\end{equation}
Combining \eqref{3.6}--\eqref{3.9}, we obtain
\begin{eqnarray}\label{3.10}
&&M^{6}\int_{X}d(x)^{2}|\nabla w|^2dx\notag\\
&\le&C_2\int_{X}d(x)^{2}|\nabla w(x)|^2dx+C_2\int_{X}d(x)^{4}|\nabla q|^2dx+C_2M^{12}\int_{X}|w|^2dx.\notag\\
\end{eqnarray}

On the other hand, we have from the first equation of \eqref{2.6} that

\begin{eqnarray}\label{3.11}
&& E^{2}\int_{|x-y|\leq d(y)/4}|\nabla q|^2dx\notag\\
&\le&C_3((at)^2+d(y)^{-2})\int_{|x-y|\leq d(y)/2}|\nabla q|^2dx\notag\\
&& +C_3(E^4+d(y)^{-4}+(at)^4)\int_{|x-y|\leq d(y)/2}|q|^2dx.
\end{eqnarray}
Now taking $E=Md(y)^{-1}$ and
multiplying $d(y)^{6}$ on both sides of \eqref{3.4}, we have
\begin{eqnarray}\label{3.12}
&&M^{2}d(y)^{4}\int_{|x-y|\leq d(y)/4}|\nabla q|^2dx\notag\\
&\le& C_3((at)^2d(y)^{2}+1)\int_{|x-y|\leq d(y)/2}d(y)^{4}|\nabla q|^2dx\notag\\
&&+C_3(M^4+1+(at)^4d(y)^{4})\int_{|x-y|\leq d(y)/2}d(y)^{2}|q|^2dx.
\end{eqnarray}
Repeating the arguments in \eqref{3.6}$\sim$\eqref{3.10}, we have that
\begin{eqnarray}\label{3.13}
&&M^{2}\int_{X}d(x)^4|\nabla q|^2dx\notag\\
&\le& C_4\int_X((at)^2d(x)^{2}+1)d(x)^{4}|\nabla q|^2dx\notag\\
&&+C_4\int_X(M^4+1+(at)^4d(x)^{4})d(x)^{2}|q|^2dx\notag\\
&\le& C_5\int_X((at)^2d(x)^{2}+1)d(x)^{4}|\nabla q|^2dx\notag\\
&&+C_5\int_X(M^4+1+(at)^4d(x)^{4})d(x)^{2}|\nabla w|^2dx.
\end{eqnarray}
Combining \eqref{3.13} and \eqref{3.10}, we obtain that if $M\ge M_0$ for some $M_0>1$ then
\begin{eqnarray}\label{3.14}
&&M^{4}\int_{X}d(x)^{2}|\nabla w|^2dx+M^2\int_{X}d(x)^{4}|\nabla q|^2dx\notag\\
&\le&C_6\int_{X}((at)^4d(x)^{4})d(x)^{2}|\nabla w(x)|^2dx\notag\\
&&+C_6M^{12}\int_{X}|w|^2dx+C_6\int_X((at)^2d(x)^{2})d(x)^{4}|\nabla q|^2dx.\notag\\
\end{eqnarray}
Note that $B_{a_2}\subset\Omega_t$ and therefore
$$d(x)<\frac 1a-\frac{1}{at}<1.$$
Taking $M=(C_6+1)at$, one can eliminate $\int_{X}d(x)^{4}|\nabla q|^2dx$ and $\int_{X}d(x)^{2}|\nabla w(x)|^2dx$
on the right hand side of \eqref{3.14}. Finally, we get that
\begin{eqnarray}\label{3.15}
&&(at)^4\int_{X}d(x)^{2}|\nabla w|^2dx+\int_{X}d(x)^{4}|\nabla q|^2dx\notag\\
&\le&C_{7}(at )^{12}\int_{X} |w|^2 dx.
\end{eqnarray}
It is no harm to add $\int_{X}d(x)^{2}|q|^2dx$ to the right hand side of \eqref{3.15} since $\int_Xd(x)^2|q|^2dx\le\int_Xd(x)^2|\nabla w|^2dx$. We then obtain \eqref{3.1}.
\eproof

\section{Proof of Theorem \ref{thm1.1}}\label{sec4}
\setcounter{equation}{0}

This section is devoted to the proof of the main theorem, Theorem \ref{thm1.1}.
Since $(w,p)\in (H^1(\Omega_t))^{n+1}$, the regularity theorem implies $w\in
H^2_{loc}(\Omega_t)$. Therefore, to use estimate
\eqref{2.7}, we simply cut-off $w$. So let $\chi(x)\in C^{\infty}_0
({\mathbb R}^n)$ satisfy $0\le\chi(x)\leq 1$ and
\begin{equation*}
\chi (x)=
\begin{cases}
\begin{array}{l}
0,\quad |x|\leq \frac{1}{4at},\\
1,\quad \frac{1}{2at}<|x|<\frac{1}{a}-\frac{3}{at},\\
0,\quad |x|\geq \frac{1}{a}-\frac{2}{at}.
\end{array}
\end{cases}
\end{equation*}
It is easy to see that for any multiindex
$\alpha$
\begin{equation}\label{4.1}
\begin{cases}
|D^{\alpha}\chi|=O((at)^{|\alpha|})\quad \text{if}\quad \frac{1}{4at}\le |x|\le \frac{1}{2at},\\
|D^{\alpha}\chi|=O((at)^{|\alpha|})\quad \text{if}\quad \frac{1}{a}-\frac{3}{at}\le |x|\le \frac{1}{a}-\frac{2}{at}.
\end{cases}
\end{equation}
If we choose $s<8r_1$, then $\supp(\chi)\subset B_{r_1}$, where $r_1$ is defined in Lemma \ref{lem2.1}. Therefore, applying \eqref{2.8} to $\chi w$ gives
\begin{eqnarray}\label{4.2}
&&\int (\log|x|)^{-2}\varphi^2_\beta|x|^{-n}(\beta|x|^{2}|\nabla (\chi w)|^2+\beta^3|\chi w|^2)dx\notag\\
&\leq&C\int\varphi^2_\beta |x|^{-n}|x|^{4}|\Delta (\chi w)|^2 dx.
\end{eqnarray}
Here and after, $C$ and $\tilde C$ denote general constants whose value may vary from line to line. The dependence of $C$ and $\tilde C$ will be specified whenever necessary. Next applying \eqref{2.7} to $v=\chi q$ yields that
\begin{eqnarray}\label{4.3}
&&\int\varphi^2_\beta (\log|x|)^2(|x|^{4-n}\beta|\nabla (\chi q)|^2+|x|^{2-n}\beta^3|\chi q|^2)dx\notag\\
&\leq& C\int \varphi^2_\beta(\log|x|)^{4}|x|^{6-n}|\Delta(\chi q)|^2dx.
\end{eqnarray}

Combining \eqref{4.2} and  \eqref{4.3}, we obtain that
\begin{eqnarray}\label{4.4}
&&\int_{W}(\log|x|)^{-2}\varphi^2_\beta |x|^{-n}(\beta|x|^{2}|\nabla w|^2+\beta^3|w|^2)dx\notag\\
&&+\int_{W}(\log|x|)^2\varphi^2_\beta |x|^{-n}(\beta|x|^{4}|\nabla q|^2+|x|^{2}\beta^3|q|^2)dx\notag\\
&\leq& \int\varphi^2_\beta (\log|x|)^{-2}|x|^{-n}(\beta|x|^{2}\nabla (\chi w)|^2+\beta^3|\chi w|^2)dx\notag\\
&&+ \int(\log|x|)^2\varphi^2_\beta |x|^{-n}(\beta|x|^{4}|\nabla (\chi q)|^2+\beta^3|x|^{2}|\chi q|^2)dx\notag\\
&\leq &C\int\varphi^2_\beta |x|^{-n}|x|^{4}|\Delta (\chi w)|^2 dx\notag\\
&&+C\int \varphi^2_\beta(\log|x|)^{4}|x|^{6-n}|\Delta(\chi q)|^2dx,
\end{eqnarray}
where $W=\{x:\frac{1}{2at}<|x|<\frac{1}{a}-\frac{3}{at}\}$.
Define $Y=\{x:\frac{s}{32t}=\frac{1}{4at}\le |x|\le \frac{1}{2at}=\frac{s}{16t}\}$ and $Z=\{x:\frac{1}{a}-\frac{3}{at}\le |x|\le \frac{1}{a}-\frac{2}{at}\}$.
By \eqref{2.6} and estimates
\eqref{4.1}, we deduce from \eqref{4.4} that
\begin{eqnarray}\label{4.5}
&& \int_{W}(\log|x|)^{-2}\varphi^2_\beta |x|^{-n}(\beta|x|^{2}|\nabla w|^2+\beta^3|w|^2)dx\notag\\
&&+\int_{W}(\log|x|)^2\varphi^2_\beta |x|^{-n}(\beta|x|^{4}|\nabla q|^2+|x|^{2}\beta^3|q|^2)dx\notag\\
&\leq& C\int_{W}\varphi^2_\beta |x|^{-n}|x|^{4}|\nabla q|^2 dx\notag\\
&&+C\int_{W}(\log|x|)^{4}\varphi^2_\beta |x|^{-n}|x|^{6}((at)^4|q|^2+(at)^2|\nabla q|^2)dx\notag\\
&&+C(at)^4\int_{Y\cup Z}\varphi^2_\beta |x|^{-n}|\tilde{U}|^2 dx\notag\\
&&+C(at)^4 \int_{Y\cup Z}(\log|x|)^{4}\varphi^2_\beta |x|^{2-n}|\tilde{U}|^2 dx,
\end{eqnarray}
where $|\tilde{U}(x)|^2=|x|^{4}|\nabla q|^2+|x|^{2}|q|^2+|x|^{2}|\nabla w|^2+|w|^2$ and the positive
constant $C$ only depends on $\lambda$ and $n$.

It is easy to check that there exists
$\tilde R_1>0$, depending on $n$, such that for all $\beta>0$, both
$(\log|x|)^{-2}|x|^{-n}\varphi_{\beta}^2(|x|)$ and
$(\log|x|)^{4}|x|^{-n}\varphi_{\beta}^2(|x|)$ are decreasing
functions in $0<|x|<\tilde R_1$. So we choose a small $s<8\min\{r_1,\tilde
R_1\}$. Now letting $\beta\geq \tilde{\beta}$ with $\tilde{\beta}=C(at)^2+1$,
then the first two terms on the right hand side of \eqref{4.5} can be absorbed by the left hand
side of \eqref{4.5}.  With the choices described above, we obtain from \eqref{4.5} that

\begin{eqnarray}\label{4.6}
&&\beta^3(b_1)^{-n}(\log b_1)^{-2}\varphi^2_\beta(b_1)\int_{\frac{1}{at}<|x|<b_1}|w|^2dx\notag\\
&\leq &\beta^3\int_{W}(\log|x|)^{-2}\varphi^2_\beta |x|^{-n}|w|^2dx\notag\\
&\leq &C(at)^4 \int_{Y\cup Z}(\log|x|)^{4}\varphi^2_\beta |x|^{-n}|\tilde{U}|^2 dx\notag\\
&\leq &C(at)^4(\log b_2)^{4}b_2^{-n}\varphi^2_\beta(b_2)\int_{Y}|\tilde{U}|^2 dx\notag\\
&&+C(at)^4(\log b_3)^{4}b_3^{-n}\varphi^2_\beta(b_3)\int_{Z}|\tilde{U}|^2 dx,
\end{eqnarray}
where $b_1=\frac{1}{a}-\frac{8}{at}$, $b_2=\frac{1}{4at}$ and $b_3=\frac{1}{a}-\frac{3}{at}$.

Using \eqref{3.1}, we can control the $|\tilde{U}|^2$ terms on the right
hand side of \eqref{4.5}. Indeed, let $X=Y_1:=\{x:\frac{1}{8at}\le |x|\le \frac{1}{at}\}$, then we can see that
$$d(x)\ge C|x|\quad\text{for all}\quad x\in Y,$$ where $C$ an absolute constant. Therefore, \eqref{3.1} implies
\begin{eqnarray}\label{4.01}
&&\int_{Y}\left(|x|^{2}|\nabla w|^2+|x|^{4}|\nabla q|^2+|x|^{2}|q|^2\right)dx\notag\\
&\le&C\int_{Y_1}\left(d(x)^{2}|\nabla w|^2+d(x)^{4}|\nabla q|^2+d(x)^{2}|q|^2\right)dx\notag\\
&\le&C(at)^{12}\int_{Y_1}|w|^2 dx.
\end{eqnarray}
On the other hand, let $X=Z_1:=\{x:\frac{1}{2a}\le |x|\le \frac{1}{a}-\frac{1}{at}\}$, then
$$d(x)\ge Ct^{-1}|x|\quad\text{for all}\quad x\in Z,$$ where $C$ another absolute constant. Thus, it follows from \eqref{3.1} that
\begin{eqnarray}\label{4.02}
&&\int_{Z}\left(|x|^{2}|\nabla w|^2+|x|^{4}|\nabla q|^2+|x|^{2}|q|^2\right)dx\notag\\
&\le&C(at)^4\int_{Z_1}\left(d(x)^{2}|\nabla w|^2+d(x)^{4}|\nabla q|^2dx+d(x)^{2}|q|^2\right)dx\notag\\
&\le&C(at)^{16}\int_{Z_1}|w|^2 dx.
\end{eqnarray}
Combining \eqref{4.6}, \eqref{4.01}, and \eqref{4.02} implies that
\begin{eqnarray}\label{4.7}
&&b_1^{-2\beta-n}(\log b_1)^{-4\beta-2}\int_{\frac{1}{at}<|x|<b_1}|w|^2dx\notag\\
&\leq & C(at)^{16}b_2^{-2\beta-n}(\log b_2)^{-4\beta+4}\int_{Y_1}|w|^2 dx\notag\\
&&+C(at)^{20}b_3^{-2\beta-n}(\log b_3)^{-4\beta+4}\int_{Z_1}|w|^2 dx.
\end{eqnarray}

Replacing $2\beta+n$ by $\beta$, \eqref{4.7} becomes
\begin{eqnarray}\label{4.8}
&&b_1^{-\beta}(\log b_1)^{-2\beta+2n-2}\int_{\frac{1}{at}<|x|<b_1}|w|^2dx\notag\\
&\leq & C(at)^{16}b_2^{-\beta}(\log b_2)^{-2\beta+2n+4}\int_{Y_1}|w|^2 dx\notag\\
&&+C(at)^{20}b_3^{-\beta}(\log b_3)^{-2\beta+2n+4}\int_{Z_1}|w|^2 dx.
\end{eqnarray}
Dividing $b_1^{-\beta}(\log b_1)^{-2\beta+2n-2}$ on the both sides
of \eqref{4.8} and noting that $\beta\geq n+2>n-1$, i.e., $2\beta-2n+2>0$, we get
\begin{eqnarray}\label{4.9}
&&\int_{|x+\frac{b_4x_0}{t}|<\frac{1}{at}}|w(x)|^2dx\notag\\
&\leq &\int_{\frac{1}{at}<|x|<b_1}|w(x)|^2dx\notag\\
&\leq & C(at)^{16}(\log (4at))^{6}(b_1/b_2)^{\beta}\int_{Y_1}|w|^2dx\notag\\
&&+C(at)^{20}(b_1/b_3)^\beta [\log b_1/\log b_3]^{2\beta-2n+4}\int_{Z_1}|w|^2 dx\notag\\
&\leq & C(at)^{16}(\log (4at))^{6}(4t)^{\beta}\int_{|x|<\frac{1}{at}}|w(x)|^2dx+C(at)^{20}(b_1/b_5)^{\beta}\int_{Z_1}|w(x)|^2 dx,\notag\\
\end{eqnarray}
where $b_4=\frac{1}{a}-\frac{10}{at}$ and $b_5=\frac{1}{a}-\frac{6}{at}$.
In deriving the third inequality above, we use the fact that
$$
(\frac{b_5}{b_3})(\frac{\log b_1}{\log b_3})^2\le 1
$$
for all $t\ge t'_0$ and $s\le\tilde R_2$, where $t'_0$ and $\tilde R_2$ are absolute constants. So we pick $s<\min\{8r_1,8\tilde R_1,\tilde R_2\}$ and fix it from now on. We observe that $s$ depends only on $n$.

From \eqref{4.9}, \eqref{2.2} and the definition of $w(x)$, the change of variables $y=atx+x_0$ leads to
\begin{eqnarray}\label{4.10}
M(10)&\leq&  Ct^{16}(\log (4at))^{6}(4t)^{\beta}\int_{|y-x_0|< 1}|u(y)|^2dy+C(\lambda^2\omega_n) t^{20+n}\notag\\
& \leq &  C(4t)^{\beta+22}\int_{|y-x_0|< 1}|u(y)|^2dy+Ct^{20+n}(\frac{t}{t+2})^{\beta}\notag\\
&\le&C(4t)^{2\beta}\int_{|y-x_0|< 1}|u(y)|^2dy+Ct^{20+n}(\frac{t}{t+2})^{\beta},
\end{eqnarray}
where $\omega_n$ is the volume of the unit ball and thus $C$ depends on $\lambda$, $n$.
It should be noted that \eqref{4.10} holds for all $t\ge t''_0$,
$\beta\ge\tilde\beta (\ge 22)$, where $t''_0$ depends only on $n$. For
simplicity, by denoting
\begin{equation*}
A(t)=2\log 4t,\quad B(t)=\log (\frac{t+2}{t}),
\end{equation*}
\eqref{4.10} becomes
\begin{equation}\label{4.11}
M(10)\leq C\Big{\{}\exp(\beta A(t))\int_{|y-x_0|<1}|u(y)|^2dy+t^{20+n}\exp(-\beta B(t))\Big{\}}.
\end{equation}

Now, we consider two cases. If
$$ \exp(\tilde{\beta} A(t))\int_{|y-x_0|<1}|u(y)|^2dy\geq t^{20+n}\exp(-\tilde{\beta} B(t)),$$
then we have
\begin{eqnarray}\label{4.12}
\int_{|y-x_0|< 1}|u(y)|^2dy&\geq& t^{20+n}\exp(-\tilde\beta A(t)-\tilde\beta B(t))\notag\\
&\ge& t^{20+n}\left(\frac{t+2}{t}\right)^{-\tilde\beta}(4t)^{-2\tilde\beta}\notag\\
&\ge& (4t)^{-3\tilde{\beta}}\ge\exp(-Ct^2\log t),
\end{eqnarray}
where $C$ depends on $\lambda$ and $n$ and $t\ge t'''_0$.

On the other hand, if
$$\exp(\tilde{\beta} A(t))\int_{|y-x_0|<1}|u(y)|^2dy\leq t^{20+n}\exp(-\tilde{\beta} B(t)),$$
then we can pick a $\beta>\tilde\beta$ such that
$$
\exp(\beta A(t))\int_{|y-x_0|< 1}|u(y)|^2dy=t^{20+n}\exp(-\beta B(t)).
$$
Using such $\beta$, we obtain from \eqref{4.11} that
\begin{eqnarray}\label{4.13}
M(10)&\leq& C\exp(\beta A(t))\int_{|y-x_0|< 1}|u(y)|^2dy\notag\\
&=& C\left(\int_{|y-x_0|<1}|u(y)|^2dy\right)^{\tau}(t^{20+n})^{1-\tau},
\end{eqnarray}
where $\tau=\frac{B(t)}{A(t)+B(t)}$.
Thus, \eqref{4.13} implies that
\begin{equation}\label{4.14}
t^{20+n}\leq \left(\int_{|y-x_0|< 1}|u(y)|^2dy\right)\left(\frac{t^{20+n}C}{M(10)}\right)^{1/\tau}.
\end{equation}
In view of the formula for $\tau$, we can see that
$$
\frac{1}{\tau}\log\left(\frac{t^{20+n}C}{M(10)}\right)=\frac{2\log(4t)+\log(1+(2/t))}{\log(1+(2/t))}\log\left(\frac{t^{20+n}C}{M(10)}\right)\le\tilde\beta
$$
for all $t>\tilde t$. It suffices to choose $\tilde t\ge\max\{t'_0,t''_0,t'''_0\}$. It is obvious that $\tilde t$ depends on $\lambda$, $n$, and $M(10)$. Therefore, we get from \eqref{4.14} that
\begin{equation}\label{4.15}
\int_{|y-x_0|< 1}|u(y)|^2dy\ge t^{20+n}\exp(-Ct^2),
\end{equation}
where $C$ depends on $\lambda$ and $n$. Theorem~\ref{thm1.1} now follows from \eqref{4.12} and \eqref{4.15}.\eproof

\section*{Acknowledgements}
Lin and Wang were supported in part by the National Science Council
of Taiwan. Uhlmann was partially supported by NSF, a Chancellor
Professorship at UC Berkeley and a Senior Clay Scholarship at MSRI.

\end{document}